\documentclass[a4paper,preprint,showkeys]{amsproc}
\usepackage{mathrsfs}
\usepackage{amsmath,booktabs}
\usepackage[T4, OT1]{fontenc}
\usepackage[arrow, matrix, curve]{xy}
\usepackage{newunicodechar}
\usepackage{graphicx}
\usepackage{bm}
\usepackage[dvips]{epsfig}
\usepackage{rotating}
\usepackage{amsmath,amssymb,amsfonts,amsthm}
\usepackage{amsmath,amscd}
\usepackage{subfig}
\usepackage{caption}
\usepackage{pst-all}
\usepackage{dcolumn}
\usepackage{amsthm}
\usepackage{dsfont}
\usepackage{amssymb}
\usepackage[hyphens]{url} 
\usepackage[dvips]{graphicx}
\usepackage{mathrsfs}
\makeatletter
\@namedef{subjclassname@2020}{\textup{2020} Mathematics Subject Classification}
\makeatother
\theoremstyle{plain}
 \newtheorem{thm}{Theorem}[section]
 
 \newtheorem{lem}{Lemma}[section]
 
\theoremstyle{definition}
 \newtheorem{exm}{Example}[section]
 
\theoremstyle{remark}
 \newtheorem{rem}{Remark}[section]
 \numberwithin{equation}{section}

\renewcommand{\leq}{\leqslant}
\newcommand{\eg}{e.\,g.\ }
\newcommand{\ie}{i.\,e.\ }
\renewcommand{\setminus}{\smallsetminus}
\title[Volume of Hyperspheres in finite and infinite dimensions]{A novel method of computing the Volume of \\Hyperspheres in finite and infinite dimensions.}

\author[Belardinelli]{Cyril Belardinelli} 
\address{ 
Kantonsschule Solothurn\\
Solothurn\\
Switzerland}
\email{cyril.belardinelli@ksso.ch}
\begin{document}
\vspace{18mm} \setcounter{page}{1} \thispagestyle{empty}
\begin{abstract}
In the present article, the volume of the hypersphere in n-dimensional euclidean space is recalculated in a rather original way by using the theory of generalized functions (tempered distributions).\\
The calculation is performed by applying the integral representation of the Heaviside unit step function and furthermore by using the distributional Fourier transform of general power-law functions. The same method (added by a regularization procedure of the occurring infinite-dimensional integral via the Riemann-zeta function) is applied in the case of infinite dimensions giving rise to some sort of anti-measure of spheres in the Hilbert space of square-summable complex sequences.
\end{abstract}
\maketitle
\section{\label{histo}Historical Remarks}
In this article, we recompute the volume of the n-ball in a direct way using generalized functions. Obviously, the interest lies in the methodology and not in the formula itself, as it has been known for a long time. The formula for the volume is used frequently in many branches of physics and mathematics, ranging from Quantum Field Theory (\eg when applied to dimensional regularization, see \eg \cite{peskin:1995}) to Probability Theory. The volume formula has been derived for the first time by the relatively unknown Swiss mathematician L. Schl\"afli in the mid 19th century (1850-52). The formula is exposed in Schl\"afli's opus magnum \emph{Theorie der vielfachen Kontinuit\"at} \cite{graf:1901}. This work has been published posthumously only in 1901 after it has been rejected for publication several times. However, in the years between 1850-1900 most of the basic results of n-dimensional geometry have already been introduced and generally accepted. It is mainly due to this publication delay that Schl\"afli's results went largely unnoticed.
In fact, there are no comments at all on his results  
in the classical historical literature such as \eg \cite{cajori:1919, boyer:1968, kline:1973}.\\
Nevertheless, Schl\"afli must be considered (together with H. Grassmann, A. Cayley, B. Riemann and to a certain extent also W. R. Hamilton) as one of the co-founders of n-dimensional geometry. The main results have been developped by these authors \cite{cayley:1846, grassmann:1844} between the years 1844-1854.\\
\section{\label{intro}Preliminaries}
\subsection{Fourier transform of power-law functions}
The present work makes use of the theory of generalized functions (also known as tempered distributions). More precisely, this paper is basically concerned with an application of a recently derived formula for the distributional Fourier transform of a general power-law function $f(x)=x^{\alpha}$, (with a real exponent $\alpha$). This formula has been derived in \cite{belardinelli:2021} for the purpose of defining derivatives of fractional order. The formula for real (\ie not integer) exponents seems not to be known in the literature. In any case it does not appear in classical works about generalized functions such as \eg \cite{schwartz:1950, schwartz:1951, gelfand:1964, vladimirov:2002, mikusinski:1973}.\\
The classical Fourier transform of the function $f(x)=x^{\alpha}$ exists only for those exponents $\alpha$ which lie in the interval between minus one and zero. By contrast, the Fourier transform of the corresponding tempered distribution $\mathrm{Pf}\,x^{\alpha}$ is well-defined for every real exponent $\alpha$. Note that the abbreviation \emph{Pf} denotes Hadamard's Partie fini. 
The precise definition of $\mathrm{Pf}\,x^{\alpha}$ is given by: 
\begin{equation}
\mathrm{Pf}\,x^{\alpha}:=\mathrm{finite\,part\,of}\lim_{\epsilon\to0^{+}}\int\limits_{|x|>\epsilon}\varphi(x)\,x^{\alpha}\mathrm{d}x
\end{equation}
The formula for $\mathscr{F}\,\mathrm{Pf}\,x^{\alpha}$ needed below is given by the following Theorem:
\begin{thm}
\label{theorem1}
Let $f(x)=x^{\alpha}$, ($\alpha \in \mathbb{R}$) be the monomial which generates the regular distribution $\mathrm{Pf} \,x^{\alpha}$. We have then:\\
\begin{equation}
\mathscr{F}\,\mathrm{Pf}\,x^{\alpha}=\begin{cases}
2\pi \,i^{\alpha}\,\delta^{(\alpha)}(x)&\alpha\in\{0,1,2,3...\}\\
\frac{\pi\,i^{\alpha}}{(-\alpha-1)!}\,\,\mathrm{Pf}\,x^{-\alpha-1}\,\mathrm{sgn}{(x)}&\alpha\in\{-1,-2,-3...\}\\
T[\varphi]&\alpha\in\mathbb{R}\setminus\mathbb{Z}
\end{cases}
\end{equation}
where $\delta^{(\alpha)}(x)$ denotes the $\alpha^{\mathrm{th}}$ derivative of the Dirac-distribution $\delta(x)$.
$T$ denotes the regular distribution given by:
\begin{equation*}
T=\Gamma(\alpha+1)\,\left(\mathrm{e}^{-i\pi\alpha}\,\mathrm{Pf}\,\frac{1}{(-ix)^{\alpha+1}}+\mathrm{Pf}\,\frac{1}{(ix)^{\alpha+1}}\right)
\end{equation*}
\end{thm}
\begin{proof}
See \cite{belardinelli:2021}
\end{proof}
In the specific case of integer and half integer exponents Theorem \ref{theorem1} reduces to:
\begin{equation}
\label{f0}
\mathscr{F}\,\mathrm{Pf}\,x^{-m}=\mathrm{Pf}\,\begin{cases}
\frac{\pi}{i^{m}(m-1)!}\,x^{m-1}\,\mathrm{sgn}{(x)}& m\in\{1,2,3...\}\\
&\\
\sqrt{\pi}\,\,2^{m+\frac{1}{2}}\,(-2m)!!\,\Theta(-x)\,(-i|x|)^{m-1}& m\in\{\pm \frac{1}{2},\pm \frac{3}{2},\pm \frac{5}{2},...\}
\end{cases}
\end{equation}
Where $\Theta(\cdot)$ denotes Heaviside's unit step function defined by:
\begin{equation*}
\Theta(x)=\begin{cases}
1&\text{if} \quad x>0\\
\frac{1}{2}&\text{if} \quad x=0\\
0&\text{if} \quad x<0
\end{cases}
\end{equation*}
The double factorial in \eqref{f0} is defined in the usual manner:
\begin{equation*}
n!!=n\cdot(n-2)\cdot(n-4)\cdots2\,\, (1\mathrm{,respectively})
\end{equation*}
Note that the double factorial is also defined for negative odd arguments. The values can be evaluated by virtue of the formula (restricted to odd numbers $n$):
\begin{equation*}
(-n)!!\,n!!=(-1)^{\frac{n-1}{2}}\,n
\end{equation*}
A few examples are given by:
\begin{exm}
\begin{equation*}
(-1)!!=1\quad(-3)!!=-1\quad(-5)!!=\frac{1}{3}
\end{equation*}
\end{exm}
\subsection{Branch of the complex Logarithm}
We stress the fact that formula of Theorem \ref{theorem1} is based on the following branch of the complex logarithm:
\begin{align*}
\log \colon \mathbb{C}^{\times}&\to
\mathbb{C}\\
r\mathrm{e}^{i\varphi} &\mapsto \log{r}+i\varphi
\end{align*}
where $\mathbb{C}^{\times}:=\{r\mathrm{e}^{i\varphi}|\,r>0,\,-\pi\leq\varphi<\pi \}$\\
\begin{rem}
As a consequence of the chosen branch we have the following identities which we use in the paper.
\begin{equation}
\label{branch_log}
\log{i}=i\frac{\pi}{2}\mathrm{,} \quad \log{(-i)}=-i\frac{\pi}{2}\mathrm{,} \quad \log{(-1)}=-i\pi
\end{equation}
\end{rem}
\section{\label{first_part}Volume of the n-Ball}
\label{sec.iii}
In this paper, we are concerned with the volume of n-dimensional balls in $\mathbb{R}^{n}$ defined by:
\begin{equation}
\mathbb{B}^{n}(r)=\{(x_{1},x_{2},\dots,x_n)\in \mathbb{R}^n\,|\,x_{1}^{2}+x_{2}^2+\dots+x_{n}^2\leq r^2\}
\end{equation}
The Volume of $\mathbb{B}^{n}(r)$ can be expressed in the form:
\begin{equation}
\label{f1}
\mathrm{Vol}\,(\mathbb{B}^{n}(r))=\int\limits_{\mathbb{R}^n}\mathrm{d}x_{1}\mathrm{d}x_{2}\cdots \mathrm{d}x_{n}\,\Theta(r^2-x_{1}^2-x_{2}^2-\cdots-x_{n}^2)
\end{equation}
In the next step we use the integral representation of $\Theta(x)$:
\begin{equation}
\Theta(x)=\frac{1}{2\pi i}\lim_{\epsilon\to 0}\int\limits_{\mathbb{R}}\mathrm{d}k\,\frac{\mathrm{e}^{ikx}}{k-i\epsilon}
\end{equation}
By inserting this expression into \eqref{f1} we get:
\begin{equation}
\label{f50}
\mathrm{Vol}\,(\mathbb{B}^{n}(1))=\frac{1}{2\pi i}\lim_{\epsilon\to 0}\int\limits_{\mathbb{R}^n}\mathrm{d}x_{1}\mathrm{d}x_{2}\cdots \mathrm{d}x_{n}\int\limits_{\mathbb{R}}\mathrm{d}k\,\frac{\mathrm{e}^{-k(\epsilon-i)}}{k-i\epsilon}\mathrm{e}^{-(\epsilon+ik)(x_{1}^2+x_{2}^2+\cdots+x_{n}^2)}
\end{equation}
Note that we inserted factors $\mathrm{e}^{-\epsilon(\cdot)}$ into the integral in \eqref{f50} in order to ensure convergence.
\\
According to Fubini's Theorem we obtain (leaving out the factor $\mathrm{e}^{-\epsilon k}$):
\begin{equation}
\label{f2}
\mathrm{Vol}\,(\mathbb{B}^{n}(1))=\frac{\pi^{\frac{n}{2}-1}}{2}\lim_{\epsilon\to 0}\int\limits_{\mathbb{R}}\mathrm{d}k\,\frac{\mathrm{e}^{ik}}{(\epsilon+ik)^{\frac{n}{2}+1}}
\end{equation}
Next, we reinterpret the integral in \ref{f2} as the Fourier backtransform (evaluated at $x=1$) of the following \emph{tempered distribution}:
\begin{equation}
T_{\epsilon}=\mathrm{Pf}\frac{1}{(\epsilon+ik)^{\frac{n}{2}+1}}
\end{equation}
The volume can then be rewritten in the form:
\begin{equation}
\mathrm{Vol}(\mathbb{B}^{n}(1))=\pi^{\frac{n}{2}}\,\mathscr{S}^{\prime}-\lim_{\epsilon\to 0}\left[\mathscr{F}^{-1}T_{\epsilon}\,(x=1)\right]
\end{equation}
Where the notation $\mathscr{S}^{\prime}-\lim_{\epsilon\to 0}$ denotes the limit within the space of distributions $\mathscr{S}^{\prime}$.
Due to the continuity of the operator $\mathscr{F}^{-1}$ we are allowed to write:
\begin{equation}
\label{f3}
\mathscr{S}^{\prime}-\lim_{\epsilon\to 0}\,[\mathscr{F}^{-1}T_{\epsilon}]=\mathscr{F}^{-1}[\mathscr{S}^{\prime}-\lim_{\epsilon\to 0}T_{\epsilon}]
\end{equation}
The limit on the right hand side of \eqref{f3} can be determined by using the following Lemma \cite{belardinelli:2021}:
\begin{lem}
\label{lemma1}
\begin{equation*}
\mathscr{S}^{\prime}-\lim_{\epsilon\to0^{+}}\frac{1}{(\epsilon+ix)^{m+1}}=
\begin{cases}
\pi\,\frac{i^{m}}{m!}\,\delta^{(m)}(x)+\mathrm{Pf}\,(ix)^{-m-1}&\,m\in\{1,2,3...\}\\
&\\
\mathrm{Pf}\,(ix)^{-m-1}&\,m\in\{\pm \frac{1}{2},\pm \frac{3}{2},\pm \frac{5}{2},...\}\\
\end{cases}
\end{equation*}
\end{lem}
\begin{proof}
See \cite{belardinelli:2021}
\end{proof}
Applying formula \eqref{f0} we obtain:
\begin{equation*}
\mathscr{F}^{-1}\left[\mathscr{S}^{\prime}-\lim_{\epsilon\to 0}T_{\epsilon}\right]=
\begin{cases}
\frac{1}{m!}&\,m\in\{1,2,3,...\}\\
&\\
\frac{2^{m+\frac{1}{2}}i^{-m-1}(-i)^m}{\sqrt{\pi}}(-2m-2)!!&\,m\in\{\frac{1}{2},\frac{3}{2},\frac{5}{2},...\}\\
\end{cases}
\end{equation*}
\begin{rem}
The expression for half integer values of $m$ can be simplified by using the following identities (valid for $m\in\left\{\frac{1}{2},\frac{3}{2},\frac{5}{2},...\right\}$):
\begin{equation*}
(-2m-2)!!=\frac{(-1)^{m+\frac{1}{2}}}{(2m)!!}
\end{equation*}
and 
\begin{equation*}
(-1)^{m+\frac{1}{2}}(-i)^{m}i^{-m-1}=1
\end{equation*}
\end{rem}
The simplified formula reads then:
\begin{equation*}
\mathscr{F}^{-1}\left[\mathscr{S}^{\prime}-\lim_{\epsilon\to 0}T_{\epsilon}\right]=
\begin{cases}
\frac{1}{m!}&\,m\in\{1,2,3,...\}\\
&\\
\frac{2^{m+\frac{1}{2}}}{\sqrt{\pi}(2m)!!}&\,m\in\{\frac{1}{2},\frac{3}{2},\frac{5}{2},...\}\\
\end{cases}
\end{equation*}
\\
The final result for the volume formulated in terms of space dimension $n$ is then:
\begin{equation*}
\mathrm{Vol}\,(\mathbb{B}^{n}(1))=
\begin{cases}
\frac{\pi^{\frac{n}{2}}}{(\frac{n}{2})!}&\,n\in\{2,4,...\}\\
&\\
2\frac{(2\pi)^{\frac{n-1}{2}}}{(n)!!}&\,n\in\{1,3,5,...\}\\
\end{cases}
\end{equation*}
which can be subsumed under the well-known formula:
\begin{equation*}
\mathrm{Vol}\,(\mathbb{B}^{n}(1))=\frac{\pi^{\frac{n}{2}}}{\Gamma(\frac{n}{2}+1)}\mathrm{,}\quad n\in\{1,2,3...\}
\end{equation*}
\section{Volume of Spheres in infinite dimensions}
In the present section we try to calculate the volume of balls in the Hilbert space  of square-summable complex sequences $\ell^{2}(\mathbb{C}):=\{(c_n):c_n \in \mathbb{C} \,|\sum_{n=1}^{\infty}|c_{n}|^2\leq \infty\}$. However, it is known that a $\sigma$-additive Lebesgue-type measure which is translationally and rotationally invariant does not exist in a Hilbert space of infinite dimensions (See \eg \cite{mazzucchi:09}). The volume we are calculating in this section does therefore not represent a measure in the usual sense $\grave{a}$ la Lebesgue.\\
The volume to be calculated is formally given by an infinite-dimensional integral which we regularize by Riemann's $\zeta$-function. We remark that an analogous calculation is given in \cite{belardinelli:2019} where one calculates a certain path integral $\grave{a}$ la Feynman by using the integral representation of delta-function $\delta(x)$ and applying a $\zeta$-regularization.   
\\
\\
The ball of radius $r$ within $\ell^{2}(\mathbb{C})$ is defined by: 
\begin{equation}
\mathbb{B}_{\mathbb{C}}^{\infty}(r):=\{(c_{1},c_{2},\dots)\in\ell^{2}(\mathbb{C}) \,|\sum_{n=1}^{\infty}|c_{n}|^2\leq r^2\}
\end{equation}
Analogous to \eqref{f1} we have:
\begin{equation}
\mathrm{Vol}\,(\mathbb{B}_{\mathbb{C}}^{\infty}(r))=\int\limits_{\mathbb{C}\times\mathbb{C}\times\cdots}\prod_{n=1}^{\infty}\mathrm{d}c_{n}\,\Theta(r^2-|c_{1}|^2-|c_{2}|^2-\cdots)
\end{equation}
Where we define the integration over the complex variable $c_n=x_{n}+iy_n$ in the following canonical way:
\begin{equation} 
\int_{\mathbb{C}} \mathrm{d}c_{n}:=\int_{\mathbb{R}^2} \mathrm{d}x_{n} \mathrm{d}y_{n}
\end{equation}
Next, we proceed along the same track as in Sec. \ref{sec.iii}
\begin{equation*}
\mathrm{Vol}\,(\mathbb{B}_{\mathbb{C}}^{\infty}(r))=\frac{1}{2\pi i}\lim_{\epsilon\to 0}\int\limits_{\mathbb{C}\times\mathbb{C}\times\cdots}(\mathrm{d}c_{n})\int\limits_{\mathbb{R}}\mathrm{d}k\,\frac{\mathrm{e}^{ikr^2}}{k-i\epsilon}\,\mathrm{e}^{-(\epsilon+ik)(|c_{1}|^2+|c_{2}|^2+\cdots}
\end{equation*}
Where we use the short notation $(\mathrm{d}c_{n}):=\prod_{n=1}^{\infty}\mathrm{d}c_{n}$\\
\\
By the same procedure as in Sec. \ref{sec.iii} we interchange the order of integration to obtain:
\begin{equation}
\label{f10}
\mathrm{Vol}\,(\mathbb{B}_{\mathbb{C}}^{\infty}(r))=\frac{1}{2\pi i}\lim_{\epsilon\to 0}\int\limits_{\mathbb{R}}\mathrm{d}k\,\frac{\mathrm{e}^{ikr^2}}{k-i\epsilon}\int\limits_{\mathbb{C}\times\mathbb{C}\times\cdots}(\mathrm{d}c_{n})\,\mathrm{e}^{-(\epsilon+ik)(|c_{1}|^2+|c_{2}|^2+\cdots}
\end{equation}
The latter integral can be evaluated formally:
\begin{equation}
\label{product}
\int\limits_{\mathbb{C}\times\mathbb{C}\times\cdots}(\mathrm{d}c_{n})\,\mathrm{e}^{-(\epsilon+ik)(|c_{1}|^2+|c_{2}|^2+\cdots}=\prod_{n=1}^{\infty}\frac{\pi}{\epsilon+ik}
\end{equation}
We regularize the latter infinite product by means of Riemann's $\zeta$-function \cite{titchmarsh:1951, ahlfors:78}. We stress the fact that the same regularization procedure has already been applied successfully in evaluating Feynman path integrals with quadratic Lagrangians. See \eg \cite{belardinelli:2019}.\\
The Logarithm of the product in Eq.~\eqref{product} reads formally :
\begin{equation*}
\log{\prod_{n=1}^{\infty}\frac{\pi}{\epsilon+ik}}=\left(\log{\frac{\pi}{\epsilon+ik}}\right)\,\sum_{n=1}^{\infty}1
\end{equation*}
By applying the $\zeta$-regularized sum (See \eg \cite{belardinelli:2019}):
\begin{equation*}
\sum_{n=1}^{\infty}1=\zeta(0)=-\frac{1}{2}
\end{equation*}
we obtain:
\begin{equation*}
\prod_{n=1}^{\infty}\frac{\pi}{\epsilon+ik}=\sqrt{\frac{\epsilon+ik}{\pi}}
\end{equation*}
Finally, we insert the latter expression into \eqref{f10} which gives:
\begin{equation*}
\mathrm{Vol}\,(\mathbb{B}_{\mathbb{C}}^{\infty}(r))=\frac{1}{2\pi^{\frac{3}{2}}}\lim_{\epsilon\to 0}\int\limits_{\mathbb{R}}\mathrm{d}k\,\frac{\mathrm{e}^{ikr^2}}{\sqrt{\epsilon+ik}}
\end{equation*}
The latter expression can be rewritten formally as:
\begin{equation*}
\mathrm{Vol}\,(\mathbb{B}_{\mathbb{C}}^{\infty}(r))=\frac{1}{\sqrt{\pi}}\,\mathscr{S}^{\prime}-\lim_{\epsilon\to 0}\left[\mathscr{F}^{-1}(\epsilon+ik)^{-\frac{1}{2}}\right]_{x=r^2}
\end{equation*}
Once again, the limit $\epsilon\to 0$ can be interchanged with the operator $\mathscr{F}^{-1}$ to obtain:
\begin{equation*}
\mathrm{Vol}\,(\mathbb{B}_{\mathbb{C}}^{\infty}(r))=\frac{1}{\sqrt{i\pi}}\,\left[\mathscr{F}^{-1}\frac{1}{\sqrt{k}}\right]_{x=r^2}
\end{equation*}
According to \eqref{f0}:
\begin{equation*}
\mathscr{F}^{-1}\frac{1}{\sqrt{k}}=\frac{1}{\sqrt{\pi}\,0!!}\,\Theta(x)\,(-i|x|)^{-1/2}
\end{equation*}
Which leads to to the result:
\begin{equation}
\label{f30}
\mathrm{Vol}\,(\mathbb{B}_{\mathbb{C}}^{\infty}(r))=\frac{1}{\pi r}
\end{equation}
\begin{rem}
Formally, the same result is obtained when one inserts the value $n=-1$ into formula:
\begin{equation*}
\mathrm{Vol}\,(\mathbb{B}^{n}(r))=\frac{\pi^{\frac{n}{2}}}{\Gamma(\frac{n}{2}+1)}\,r^n
\end{equation*}
\end{rem}
Formula \eqref{f30} states that spheres with a greater radius have less volume. Evidently, this feature is incompatible with the usual properties of a Lebesgue-measure. It must therefore represent something different than a usual measure. 
\appendix
\section{}
We prove the following elementary identity:
\begin{equation*}
(-1)^{m+\frac{1}{2}}(-i)^{m}i^{-m-1}=1\mathrm{,}\quad m\in\left\{\frac{1}{2},\frac{3}{2},\frac{5}{2},...\right\}
\end{equation*}
\begin{proof}
Rewriting the left hand side of the latter expression:
\begin{equation*}
\mathrm{e}^{(m+\frac{1}{2})\log{(-1)}}\mathrm{e}^{m\log{(-i)}}\mathrm{e}^{-(m+1)\log{i}}
\end{equation*}
According to \ref{branch_log} we get:
\begin{equation*}
\mathrm{e}^{(m+\frac{1}{2})(-i\pi)}\mathrm{e}^{m(-i\frac{\pi}{2})}\mathrm{e}^{-(m+1)i\frac{\pi}{2}}
=\mathrm{e}^{i2m\pi}\,\mathrm{e}^{-i\pi}=1
\end{equation*}
\end{proof}
\appendix*
\newpage

\bibliography{References_paper_volume_hypersphere}

\begin{thebibliography}{10}

\bibitem{ahlfors:78}
L.~V. Ahlfors.
\newblock {\em Complex analysis, An introduction to the theory of analytic
  functions of one complex variable, third ed.}
\newblock McGraw-Hill Book Co., New York, 1978.

\bibitem{belardinelli:2019}
C.~Belardinelli.
\newblock {\em arXiv:1908.05226}, 2019.

\bibitem{belardinelli:2021}
C.~Belardinelli.
\newblock {\em arXiv:2105.00276}, 2021.

\bibitem{boyer:1968}
C.~B. Boyer.
\newblock {\em A History of Mathematics.}
\newblock John Wiley \& Sons, Inc., 1968.

\bibitem{cajori:1919}
F.~Cajori.
\newblock {\em A history of Mathematics, 2nd ed.}
\newblock The Macmillan Company, London, 1919.

\bibitem{cayley:1846}
A.~Cayley.
\newblock On linear transformations.
\newblock {\em Cambridge and Dublin mathematical journal}, 1:104--122, 1846.

\bibitem{grassmann:1844}
H.~Grassmann.
\newblock {\em Lineale Ausdehnungslehre}, volume~1.
\newblock Verlag von Otto Wigand, 1844.

\bibitem{gelfand:1964}
G.~E.~Shilov I.~M.~Gel'fand.
\newblock {\em Generalized Functions, Volume 1: Properties and Operations},
  volume~1.
\newblock AMS Chelsea Publishing: An Imprint of the American Mathematical
  Society, 1964.

\bibitem{kline:1973}
M.~Kline.
\newblock {\em Mathematical Thought from Ancient to Modern Times.}
\newblock Oxford University Press, 1973.

\bibitem{mazzucchi:09}
S.~Mazzucchi.
\newblock {\em Mathematical Feynman Path Integrals And Their Applications}.
\newblock World Scientific, 2009.

\bibitem{mikusinski:1973}
R.~Sikorski P.~Antosik, J.~Mikusinski.
\newblock {\em Theory of Distributions, The Sequential Approach.}
\newblock Elsevier Scientific Publishing Company, Amsterdam, 1973.

\bibitem{peskin:1995}
M.~E. Peskin and D.~V. Schr\"oder.
\newblock {\em Am Introduction to Quantum Field Theory.}
\newblock Perseus Books Publishing L. L. C., 1995.

\bibitem{graf:1901}
L.~Schl\"afli.
\newblock {\em Theorie der vielfachen Kontinuit\"at.}
\newblock Herausgegeben im Auftrag der Schweiz. Naturforschenden Gesellschaft
  von J. H. Graf, Bern, 1901.

\bibitem{schwartz:1950}
L.~Schwartz.
\newblock {\em Th\'eorie Des Distributions. Tome I}.
\newblock Herrmann Paris, 1950.

\bibitem{schwartz:1951}
L.~Schwartz.
\newblock {\em Th\'eorie Des Distributions. Tome II}.
\newblock Herrmann Paris, 1951.

\bibitem{titchmarsh:1951}
E.~C. Titchmarsh.
\newblock {\em The Theory of the Riemann Zeta-function}.
\newblock Oxford University Press, 1951.

\bibitem{vladimirov:2002}
V.S. Vladimirov.
\newblock Methods of the theory of generalized functions.
\newblock {\em Taylor \& Francis, London 2002}, pages 19--21, 2002.

\end{thebibliography}
\bibliographystyle{plain}
\end{document}